\documentclass[a4paper,12pt,oneside]{article}
\usepackage{a4wide}
\usepackage{t1enc}
\usepackage{amsfonts}
\usepackage{theorem}
\usepackage{amsmath,amssymb}
\usepackage{color}
\usepackage[dvips]{graphicx}
\usepackage{psfrag}
\usepackage{ae,aecompl}

\newcommand{\cont}[1]{\mathcal{C}^{#1}}
\newcommand{\eps}{\varepsilon}
\newcommand{\vp}{\varphi}
\newcommand{\pbe}{\displaystyle}

\newcommand{\ci}{\mathrm{i}}

\newcommand{\tor}{\mathbb{T}}
\newcommand{\rr}{\mathbb{R}}
\newcommand{\zz}{\mathbb{Z}}
\newcommand{\nn}{\mathbb{N}}

\newcommand{\cte}{\mathrm{cte} \,}

\newcommand{\ocal}{\mathcal{O}}

\newcommand{\LL}{\mathcal{L}}

\newcommand{\tet}{\widetilde{\theta}}

\newtheorem{theorem}{Theorem}[section]

\newtheorem{prop}[theorem]{Proposition}

\newtheorem{rem}[theorem]{Remark}

\begin{document}

 \title{A geometric mechanism of diffusion: Rigorous verification in a priori unstable Hamiltonian 
 systems\thanks{Supported in part by MICINN-FEDER MTM2009-06973 and CUR-DIUE 2009SGR859 grants.}}


%
\author{Amadeu Delshams$^1$ and Gemma Huguet$^{2,3}$\\
\parbox{12.5cm}{
  \small
  \begin{itemize}
  \item[$^1$]
    Dep. de Matem\`atica Aplicada I, Universitat Polit\`ecnica de Catalunya\\
    Av. Diagonal 647, 08028 Barcelona\\
    {\footnotesize \texttt{Amadeu.Delshams@upc.edu}}
  \item[$^2$]
    Centre de Recerca Matem\`{a}tica\\
    Campus de Bellaterra, Edifici C, 08193 Bellaterra (Barcelona)\\    
    {\footnotesize \texttt{Gemma.Huguet@upc.edu}}
  \item[$^3$]
     Center for Neural Science, New York University\\
     4 Washington Place, New York, NY 10003
 \end{itemize}
}}

%

\maketitle

\abstract{ 
In this paper we consider a representative \emph{a priori} unstable Hamiltonian system with $2+1/2$ degrees of freedom, to which we apply 
the geometric mechanism for diffusion introduced in the paper Delshams et al., \emph{Mem. Amer.
Math. Soc.} 2006, and generalized in  Delshams and Huguet, \emph{Nonlinearity} 2009, and provide explicit, concrete and easily verifiable conditions 
for the existence of diffusing orbits.

The simplification of the hypotheses allows us to perform explicitly
the computations along the proof, which contribute to present in an
easily understandable way the geometric mechanism of diffusion. In particular, we fully describe the construction of the scattering map and the
combination of two types of dynamics on a normally hyperbolic invariant manifold.
}


\section{Introduction}\label{sec:introduction}

The goal of this paper is to apply the geometric mechanism for
diffusion introduced in \cite{DelshamsLS06} and generalized in
\cite{DelshamsH09}, to a representative \emph{a priori}
unstable Hamiltonian system with $2+1/2$ degrees of freedom,
establishing explicit conditions for diffusion.

The phenomenon of global instability in nearly integrable Hamiltonian
systems, called \emph{Arnold diffusion}, deals essentially with the question
of what is the effect on the dynamics when an autonomous integrable
mechanical system is subject to a small periodic perturbation.
That is, whether these effects accumulate over time leading to a
large term effect or whether they average out.

For an integrable Hamiltonian system written in
action-angle variables all the trajectories lie on invariant tori,
with an associated dynamics consisting of a rigid rotation with constant frequency. For a general
perturbation of size $\eps$ of a non-degenerate integrable
Hamiltonian, the KAM Theorem (see \cite{Llave01} for a survey)
ensures the stability for most of the trajectories of the system.
More precisely, those invariant tori in the unperturbed system
$\eps=0$ having \emph{Diophantine frequencies} are preserved (they
are tori with \emph{non-resonant} frequencies), giving rise to a
Cantorian foliation of invariant tori for the perturbed system $\eps> 0$.
Thus, the set of surviving tori has a large measure but
also many gaps among them, which can be of size up to order $\sqrt{\eps}$.
However, nothing is said by the KAM theorem about the stability of
the trajectories that do not lie on the \emph{non-resonant invariant
tori}. Besides, for systems with more than two degrees of freedom
the invariant tori are not anymore an obstruction for the existence
of trajectories that may possibly drift arbitrarily far, called
\emph{diffusing orbits}.

The first description of this instability phenomenon was given by
Arnold in \cite{Arnold64} by means of his celebrated example and has
been thoroughly studied since then, using a wide range of
techniques: geometric, variational and topological (see
\cite{DelshamsGLS08} for a long list of references).

The \emph{geometric mechanism for diffusion} in
\cite{DelshamsLS06,DelshamsH09} is based on the classical Arnold
mechanism for diffusion, which consists of constructing a
\emph{transition chain}, that is, a finite sequence of whiskered
transition tori (lower dimensional invariant tori having
\emph{non-resonant frequencies} with associated stable and unstable
manifolds) with transverse heteroclinic intersections, that is, the
unstable manifold of each transition torus intersects transversally
the stable manifold of the next one.  Arnold \cite{Arnold64}
considered an integrable system with a hyperbolic component (a rotor
and a pendulum) plus a particular periodic in time perturbation,
that allowed him to construct a transition chain of \emph{primary
KAM tori}, that is, lower dimensional tori that are just a
continuation of the lower dimensional invariant tori that existed in
the unperturbed case. Using a topological argument, he proved that
diffusing orbits exist in a neighborhood of the tori in a transition chain.

Nevertheless, when one considers a generic perturbation, one faces
the problem that the gaps of size $\sqrt{\eps}$ in the foliation of
primary KAM tori are bigger than the size $\eps$ of the heteroclinic
intersection of their whiskers. Therefore, this prevents the unstable whisker of
a primary KAM torus of this foliation intersecting the stable whisker
of the next surviving primary KAM torus and one can not construct a
transition chain using only primary KAM tori. This is known in the
literature as the \emph{large gap problem} and has been solved very
recently by different methods \cite{DelshamsLS00, DelshamsLS06,DelshamsH09,%
ChengY04,ChengY09,Treschev04,PiftankinT07,GideaL06,GideaL06b}.

The strategy in \cite{DelshamsLS06, DelshamsH09} to overcome the
\emph{large gap problem} was to incorporate in the \emph{transition
chain} other invariant objects, which are not present in the
unperturbed system and are created by the resonances therein, in
order to fill the gaps between two primary KAM tori. Among them,
there are the so-called \emph{secondary KAM tori}, which are lower
invariant KAM tori topologically different from the primary ones.
The \emph{scattering map} \cite{DelshamsLS08} is the essential tool
to study the heteroclinic connections between invariant objects like
primary or secondary KAM tori.

In \cite{DelshamsLS06} it was proved the existence of Arnold
diffusion in a priori unstable Hamiltonian systems of $2+1/2$
degrees of freedom, under concrete geometric hypotheses. However, one of
them was the assumption of a non-generic condition, namely, that the
Hamiltonian was a trigonometric polynomial in the angular variables.
This latter assumption was removed in \cite{DelshamsH09} and
therefore, for these kind of Hamiltonians, the conditions required
for the geometric mechanism of diffusion were proven to be generic
in the $\cont{2}$ topology. Moreover, the removal of this hypothesis
allowed us to present the conditions for diffusion more explicitly
in terms of the original perturbation.

A strong feature of this mechanism, in contrast to other
existing ones, is that the conditions for diffusion are computable
and therefore verifiable in specific examples. Moreover, the way the
mechanism is designed allows us to give an explicit description of the
diffusing orbits. From our point of view, this fact makes this
mechanism really attracting for applications, where the computation
of the diffusing orbit in concrete systems is the cornerstone of the problem
(see \cite{DelshamsMR08}).

Although the conditions for diffusion for any concrete system are 
totally explicit and computable,
the computations needed  involve the application of several steps
 of the averaging method, the expansion in $\eps$ of a NHIM, and the
verification of the existence of non-degenerate critical points of the
Melnikov potential along the straight lines. The computation and verification
of these conditions may hide the elementary features of the geometrical method. 
Because of this, in this paper we have chosen a representative kind
of a priori unstable Hamiltonian systems, where the required hypotheses
for the application of the geometric method are trivially fulfilled.
On the contrary, the geometry of the existence of non-degenerate homoclinic
orbits to the NHIM as well as the behavior of the scattering map on the NHIM
and its interaction with the inner dynamics in the NHIM can be fully described.

%


The main result of this paper is Theorem \ref{teo:main}, which
states the existence of diffusion under very concise and easily
verifiable hypotheses, for a representative class of
priori unstable Hamiltonian systems with $2+1/2$ degrees of freedom.


Our aim is to explain the strategy of the mechanism in a clear and
understandable way, and for this reason we have chosen an
illustrative type of Hamiltonian systems which appear commonly in
the literature, for which the computations along the proof can be
performed explicitly. In Section~\ref{sec:mresult} we discuss the
reasons for the particular choices we have made. Overall, we wanted
to avoid the problem that the technical details hide the main ideas
behind the mechanism.

Another important point of this paper is that we
can give explicit expressions for the equations defining the invariant tori as well as
for the Melnikov function and the reduced Poincar\'{e} function,
which are essential for the computation of the scattering map. For those readers interested in further their
understanding of the scattering map, will find here a good illustration of its construction for a particular example,
based on geometrical considerations.

Moreover, in our description we try to present an approach that
emphasizes those points that were crucial to prove diffusion. Hence,
we describe the two different dynamics on a normally hyperbolic
invariant manifold that need to be combined to create diffusion and
we compute them explicitly.

Although this paper strongly relies on the results obtained in the
previous papers \cite{DelshamsLS06, DelshamsH09}, we have made an
effort to make it self-contained for the reader mainly interested in
the heuristic description of the mechanism and how it applies to
concrete examples. However, we accompany the exposition with precise
references to the results in \cite{DelshamsLS06, DelshamsH09} for
those readers concerned with the rigorous proofs for more general systems.


The paper is organized as follows: in Section~\ref{sec:mresult} we
introduce a representative kind of a priori unstable Hamiltonians with
$2+1/2$ degrees of freedom considered in this paper, for which we
can state our main result, Theorem~\ref{teo:main}, which establishes
the existence of a diffusing orbit for the model considered.  In
Section~\ref{sec:verif} we perform the explicit verification of the
mechanism for the Hamiltonian of Theorem~\ref{teo:main}. The
verification is structured in four parts, and includes a detailed
description of the scattering map.


\section{Set up and main result}\label{sec:mresult}

In \cite{DelshamsH09} there were given explicit conditions for the
existence of a diffusing orbit for generic \emph{a priori} unstable
Hamiltonian systems. That paper was mainly focused on proving the
genericity of the result, so although conditions were explicit, some
computational effort was required to check them. As we already
mentioned in the introduction, in this paper we plan to give sufficient
conditions, easily verifiable, which guarantee the existence of
diffusion for a general case of a priori unstable Hamiltonian
systems.

In this section, we first introduce a representative type of a priori
unstable Hamiltonian systems of $2+1/2$ degrees of freedom,
which is usually found with several variations in the literature
\cite[\textsection 7.5]{Chirikov79}\cite[\textsection 9, \textsection 12]{ChierchiaG94}%
\cite{BessiCV01,BertiB02,BertiBB03,Treschev04},
when explicit computations are performed.
We first discuss the features
and particularities of this type of systems and we finally state
Theorem \ref{teo:main}, which establishes the existence of diffusing
orbits for these systems under explicit and easily
verifiable conditions.

We consider an \emph{a priori} unstable Hamiltonian system as
introduced by Chierchia and Gallavotti \cite[sections 7.5 and 12]{ChierchiaG94},
consisting of a $2 \pi$-periodic in time perturbation of a pendulum
and a rotor described by the following non-autonomous Hamiltonian
\begin{equation} \label{eq:ham}
\begin{array}{rcl}
H_{\eps}(p,q,I,\varphi,t) & = &  \pbe H_0(p,q,I)+ \eps h(p,q,I,\varphi,t;\eps)\\
 & = & \pbe P_{\pm}(p,q) + \frac{1}{2}I^2+\eps h(p,q,I,\varphi,t;\eps).
\end{array}
\end{equation}
We notice that a motivation for the model above comes from a normal
form around a resonance of a nearly integrable Hamiltonian, and we
refer the reader to \cite{DelshamsG01,KaloshinL08} for more details.

The second term $\frac{1}{2}I^2$ in the integrable Hamiltonian
$H_0(p,q,I)$ of Hamiltonian \eqref{eq:ham} describes a \emph{rotor}
and the first one
\begin{equation}\label{eq:pend}
P_{\pm}(p,q)= \pm \left ( \frac{1}{2} p^2 + V(q) \right )
\end{equation}
a \emph{pendulum}, where the potential $V(q)$ is a $2 \pi$-periodic
function, whose non-degenerate maxima give rise to saddle points of
the pendulum \eqref{eq:pend} and therefore, to hyperbolic invariant
tori of the Hamiltonian \eqref{eq:ham} when $\eps=0$. Typically it
is assumed that the maximum of $V$ is attained at the origin $q=0$,
where $V$ is assumed to vanish, as it is case for the standard
pendulum, where
\begin{equation}\label{eq:V}
V(q)=\cos q -1.
\end{equation}
This is the simple and standard choice of potential $V(q)$ that we
are going to consider in this paper, so that
\[
P_{\pm}(p,q)=\pm\left(p^{2}/2+\cos q -1\right).
\]
The origin $(p=0,q=0)$ is a saddle point of the standard pendulum,
and its \emph{separatrix} $P_{\pm}^{-1}(0)$ for positive $p$ is
given by
\begin{equation}\label{eq:homoclinic}
q_{0}(t)= 4\arctan e^{\pm t}, \quad p_{0}(t) = 2/{\cosh t}.
\end{equation}
Notice that other choices of $V$ give rise to different separatrices
that are not usually so simple.

The negative sign in the pendulum \eqref{eq:pend} is only to
emphasize the fact that the geometric mechanism we are using does not require the
Hamiltonian $H_0$ to be positive definite, as it is the case in the
variational approach, see \cite{ChengY09}.

The term $\eps h$ in \eqref{eq:ham} is the perturbation term and
depends periodically on time and on the angular variable $\vp$, so
that $h$ can be expressed via its Fourier series in the variables
$(\vp,t)$
\begin{equation}\label{eq:fper}
h(p,q,I,\vp,t;\eps)=\pbe \sum_{(k,l) \in \zz^2}h_{k,l}(p,q,I;\eps)
e^{\ci (k \vp + l t)}.
\end{equation}


It is common in the literature to consider a perturbation
\eqref{eq:fper} depending only on the angular variables $(q,\vp,t)$,
and formed by the product of a function of the pendulum variable $q$
times a function of the angular variables $(\vp,t)$
\begin{equation}\label{eq:perturbation}
h(p,q,I,\vp,t)=f(q)g(\vp,t).
\end{equation}
This is the kind of perturbation that we are going to consider
along this paper, particularly because with this choice of $h$ the
Melnikov potential \eqref{eq:melnikov}, which will be an essential
tool for the computations along the paper, has the same harmonics as
the function $g$, and they can be computed explicitly. So we will
focus on a concrete type of Hamiltonians of the form
\begin{equation}
\label{eq:example}
H_{\eps}(p,q,I,\vp,t)=\pm\left(\frac{p^{2}}{2}+\cos q -1\right)
+\frac{I^{2}}{2} +\eps f(q)g(\vp,t),
\end{equation}
defined for any real value of $(p,q,I,\vp,t,\eps)$ and $2 \pi$-periodic in the angular variables $(q,\vp,t)$.


The function $f$ could be any $2 \pi$-periodic function. However, to easily compute the harmonics of the Melnikov potential
\eqref{eq:melnikov}, we are going to assume along this paper that
$f$ has the simple form:
\begin{equation}\label{eq:f}
f(q)=\cos q.
\end{equation}

About the choice of $f$ we would like to remark two important
features. On the one hand, thanks to the fact that $f'(0)=0$ the
problem does not require the use of the theory of normally
hyperbolic invariant manifolds to ensure the persistence of these
type of objects. This simplifies the exposition and the
computations, but since we do not assume $f(0)=0$ the problem
maintains all the richness and complexity of the \emph{large gaps
problem}.  So, although the choice \eqref{eq:f} for $f$ may seem very
restrictive, we would like to insist on the fact that the complexity
of the original problem is preserved. At the beginning of Section \ref{sec:nhim} we
discuss in detail the role of the function $f$ in the problem.

A general function $g$ is of the form
\[
g(\vp,t)=\sum_{(k,l) \in \nn^2} a_{k,l}\cos(k\vp - lt -
\sigma_{k,l}) + \tilde a_{k,l}\cos(k\vp + lt - \tilde \sigma_{k,l}),
\]
which in general has an \emph{infinite number of harmonics} in the
angles $(\vp,t)$ and where $\sigma_{k,l}, \tilde \sigma_{k,l} \in
\tor$.

Since for simplicity we will study diffusion only for $I$ positive 
along this paper, we will consider $\tilde a_{k,l}=0$, that is,
\begin{equation}\label{eq:g}
g(\vp,t)=\sum_{(k,l) \in \nn^2} a_{k,l}\cos(k\vp - lt
-\sigma_{k,l}).
\end{equation}

In a natural way, and also for simplicity, we have chosen $g$ to be
an analytic function and therefore we will assume an exponential
decay for its Fourier coefficients. That is, $|a_{k,l}| \leq e^{-
\delta |(k,l)|}$, where $\delta$ is the size of the domain of
analyticity. In this paper we simply are going to assume that they have
some general lower bound with exponential decay, that is
\[
e^{-\beta \delta |(k,l)|} \leq |a_{k,l}| \leq e^{- \delta |(k,l)|},
\]
where $1 \leq \beta < 2$. Or, equivalently, we are going to assume
\begin{equation}\label{eq:aklb1}
{\hat\alpha} {\rho}^{\beta k} {r}^{\beta l} \leq |a_{k,l}| \leq \alpha
\rho^k r^l,
\end{equation}
where $1 \leq \beta < 2$ and $0<\hat \alpha < \alpha$. Moreover, $0<\rho,r < 1$ are real numbers that will be chosen small
enough.

Although the lower bound for the coefficients $a_{k,l}$ in the above equation seems very restrictive, there are several reasons for this particular choice. 
For the more expert reader, let us say that condition \eqref{eq:aklb1} implies that big gaps of maximal size appear for all the resonances in first order with respect to 
the parameter $\eps$, without performing any step of averaging. This feature is explained thoroughly in Section \ref{sec:inner}, after equation \eqref{eq:Uex}. A second reason is that requirements
\eqref{eq:aklb1} are simple to state and verify. A generic, and, of course,
more technical, set of conditions for generic perturbations are
given explicitly in \cite{DelshamsH09}. 
When the lower bound condition for $a_{kl}$ in \eqref{eq:aklb1} are not satisfied, several steps of averaging are needed.



We can now state our main result:

\begin{theorem} \label{teo:main}
Consider a Hamiltonian of the form \eqref{eq:example}, where $f(q)$ is given by
\eqref{eq:f} and $g(\vp,t)$ is any analytic function of the form
\eqref{eq:g} with non-vanishing Fourier coefficients satisfying
\eqref{eq:aklb1}. Assume that
\begin{equation}\label{eq:clambda}
\lambda:=\left | \frac{a_{1,0}}{a_{0,1}} \right |<1/1.6 \quad\mathrm{or}\quad \lambda > 1.6,
\end{equation}
as well as $0<\rho \leq \rho^{*}$ and $0<r \leq r^{*}$, where $\rho^{*}(\lambda,\alpha, \hat \alpha, \beta)$ and $r^{*}(\lambda,\alpha, \hat \alpha, \beta)$ are small enough.

Then, for any $I_{+}^{*}>0$, there exists $\eps^{*}=\eps^{*}(I_{+}^{*})>0$ such that for any $-1/2<I_{-}<I_{+} \leq I_{+}^{*}$ and 
$0 < \eps < \eps^{*}$, there exists a trajectory
$(p(t),q(t),I(t),\vp(t))$ of the Hamiltonian \eqref{eq:ham} such
that for some $T>0$
\[ I(0) \leq I_{-}; \qquad I(T) \geq I_{+}.\]
\end{theorem}

We want to remark now that not every perturbation \eqref{eq:perturbation} gives rise to
diffusion in the action $I$. In particular, if the function
$g(\vp,t)$ in \eqref{eq:g} does not depend on $\vp$, the action $I$
is a first integral, so it does not change at all. Moreover, if
$g(\vp,t)$ does not depend on $t$, Hamiltonian \eqref{eq:example} is
autonomous and therefore $H_\eps$ is a first integral, so that only
deviations of size $\sqrt{\eps}$ are possible for the action $I$.
The same happens when the function $g(\vp,t)$ does not depend fully
on the two angular variables, but only through a integer linear
combination of them---that is, $g(\vp,t)=G(\psi)$, where $\psi=k_0
\vp - l_0 t$ is an integer combination of the angular variables
$\vp,t)$---, as can be easily checked by introducing $\psi$ as a new
angular variable. In these three cases, an infinite number of Fourier coefficients
$a_{k,l}$ of the function $g(\vp,t)$ in \eqref{eq:g} vanish. 

\begin{rem}
The condition \eqref{eq:clambda} for $a_{1,0}$ and $a_{0,1}$ jointly with the assumption $\rho^{*}, r^{*}$ being small enough, will ensure the existence of a global 
connected homoclinic manifold in Section \ref{sec:nhim}. Indeed, for $\rho^{*},r^{*}$ small enough, thanks to the upper bound \eqref{eq:aklb1}, the perturbation $g$
in \eqref{eq:g} can be approximated by its first order trigonometric polynomial 
\[ g^{[\leq 1]} (\vp,t)= a_{0,0} + a_{1,0} \cos(\vp-\sigma_{1,0}) + a_{0,1} \cos (s-\sigma_{0,1}),\]
in such a way that all the computations required for the Melnikov potential \eqref{eq:melnikov} will depend explicitly on $g^{[\leq 1]}$.

\end{rem}

\section{Proof of Theorem \ref{teo:main}}\label{sec:verif}

In the following sections, we will consider any Hamiltonian
satisfying the hypotheses of Theorem \ref{teo:main}, and
we will show how the geometric mechanism in
\cite{DelshamsLS06,DelshamsH09} can be applied to construct
diffusing orbits.

We already mentioned in the Introduction of this paper that the
geometric mechanism in \cite{DelshamsLS06,DelshamsH09} is based on
the classical Arnold mechanism for diffusion. It consists of
constructing a \emph{transition chain}, that is, a finite sequence of
whiskered transition tori such that the unstable manifold of
each torus intersects transversally the stable manifold of the next
one. As a main novelty, in \cite{DelshamsLS06, DelshamsH09} the
transition chain incorporates primary as well as secondary KAM tori
created by the \emph{resonances} (as already mentioned in the
introduction), in order to overcome the \emph{large gap problem}.

In this paper, we will try to present a description of the mechanism
that emphasizes more the geometrical aspects.
We think that this description may contribute to a
better understanding and applicability of the mechanism.

In order to prove the existence of a diffusing orbit we will
identify first a NHIM (normally hyperbolic invariant manifold)  with associated stable and unstable manifolds).
It will organize the different invariant objects involved in the
transition chain (the skeleton for the diffusing orbit).

The diffusing orbit we are looking for starts on a
point close to the NHIM and in finite time reaches another point close to
the NHIM but arbitrarily far from the original one. Of course, if
the starting point lies just on the NHIM (3D), the invariant tori
(2D) inside the NHIM act as barriers for
diffusion and the orbit is confined in a bounded domain.
Fortunately, there exists an external dynamics to the NHIM, provided
by its associated stable and unstable manifolds, which will be
essential to overcome the obstacles of the invariant tori and escape
from them, as long as the starting point does not lie on the NHIM but very
close to it. Hence, it is crucial for the mechanism that the
external dynamics does not preserve the invariant tori existing in
the NHIM. Otherwise, the orbit will be confined in a finite domain
by both the inner and the external dynamics with no possibility to
escape.

Therefore, given a Hamiltonian of the form \eqref{eq:example}, the
proof of the existence of diffusing orbits can be sketched with the
following steps: Detect the NHIM (3D) and the associated stable and
unstable manifolds (4D), determine the inner and the outer dynamics
of the NHIM as well as the invariant objects for each one, and
finally show that the outer dynamics does not preserve the invariant
objects for the inner one.

One of the novelties of this paper is the explicit description of
the outer dynamics provided by the scattering map
\cite{DelshamsLS08}. It is given by the time-$\eps$ flow of a
Hamiltonian that in first order is given by an autonomous
Hamiltonian of one degree of freedom, therefore integrable. Moreover, using a
geometric description, we are able to obtain an explicit expression
for this autonomous Hamiltonian, which is the reduced Poincar\'{e}
function \eqref{eq:rpf} with the opposite sign.

On the other hand, using averaging theory, one can show that the
Hamiltonian defining the inner dynamics can be transformed into a
normal form consisting of an integrable Hamiltonian plus a small
perturbation.

Thus, we have two dynamics defined on the NHIM that can be
approximated in suitable coordinates by one-degree of freedom
autonomous Hamiltonians. Therefore, the invariant objects are given
approximately by the level sets of these integrable Hamiltonians,
for which we provide explicit expressions.

Finally, we impose that the scattering map moves the invariant tori
for the inner dynamics, in such a way that the image under the scattering
map of each of these invariant tori intersects transversally another torus. We will show that this is a generic property.


\subsection{Part 1. Existence of a NHIM and associated stable and unstable
manifolds}\label{sec:nhim}

The first part of the proof consists of looking for a NHIM with
associated stable and unstable manifolds that intersect
transversally. In order to prove the existence of these invariant
objects, we will look for them in the unperturbed case, which is
much simpler, and then study the persistence of these objects under
the perturbation.

For $\eps=0$, Hamiltonian $H_0$ in \eqref{eq:example} consists of two
uncoupled systems: a pendulum plus a rotor. Therefore, it is clear
that the 3-dimensional manifold given by
\begin{equation}
\label{eq:nhim} \tilde \Lambda=\{(0,0,I,\vp,s): (I,\vp,s) \in \rr
\times \tor^{2}\}
\end{equation}
is an invariant manifold with associated
stable and unstable manifolds (inherited from the separatrices of
the pendulum). These manifolds coincide along a separatrix given by
\begin{equation}\label{eq:wsu}
 W^{s} \widetilde{\Lambda} = W^{u}
\widetilde{\Lambda}=\{(p_0(\tau),q_0(\tau), I, \varphi,s) : \tau \in
\mathbb{R}, I \in [-1/2,I^{*}_{+}], (\varphi,s) \in \mathbb{T}^2 \},
\end{equation}
where $(p_0(\tau),q_0(\tau))$ is the chosen orbit
\eqref{eq:homoclinic} of the pendulum,  which is homoclinic to the saddle point $p=0, q=0$.

The integrable Hamiltonian $H_0$ has a one-parameter family of
two-dimensional whiskered tori given by
\begin{equation}\label{eq:toriinv}
\mathcal{T}^{0}_{I}=\{(0,0,I,\vp,s): (\vp,s)\in\tor^{2}\},
\end{equation}
with associated frequency $(I,1)$.

When we consider the perturbation $h$, that is $\eps>0$, using
the standard theory of NHIM, see
\cite{HirschPS77,Fenichel79},
we know that for $\eps>0$ small enough, maybe restricting $\tilde
\Lambda$ to a compact subset, the manifold $\tilde \Lambda$ persists
as $\tilde \Lambda_{\eps}$, as well as its local stable and unstable
manifolds.

For any general perturbation $h$ of the form
\eqref{eq:perturbation}, if $f'(0)=0$,  the NHIM is
preserved without any deformation for any $\eps$ $\tilde
\Lambda_{\eps}=\tilde \Lambda$: $p=q=0\Rightarrow\dot p=\dot q=0$.
Moreover, if $f(0)=0$, the perturbation vanishes on $\tilde \Lambda$, so the
one-parameter family of two-dimensional invariant tori existing in
the unperturbed case remains fixed under the perturbation, as in the
Arnold's example of diffusion in \cite{Arnold64}. However, a generic
perturbation, $f(0) \neq 0$ creates gaps of size $\sqrt{\eps}$ in
the foliation of persisting invariant tori and gives rise to the \emph{large gaps problem}. See section
\ref{sec:inner} for a detailed description of the invariant objects
in this foliation.

Although it is not a generic assumption, for the clarity of exposition and
for the convenience of the reader not familiar with the theory of
NHIM, we have chosen a function $f$ in \eqref{eq:f} so that the NHIM is preserved
without deformation, that is $f'(0)=0$. Nevertheless, we want to
emphasize that this is not a necessary hypothesis for the existence
of diffusion. Indeed this was not assumed in the proof of the result
in \cite{DelshamsLS06, DelshamsH09} where NHIM theory was used. Of course, any other function $f$ satisfying the conditions
$f'(0)=0$ and $f(0) \neq 0$ will be enough for exhibiting the \emph{large gaps problem}, but we have
chosen the \emph{concrete} one \eqref{eq:f} so that the Fourier
coefficients for the Melnikov potential \eqref{eq:melnikov} can be
computed \emph{explicitly}.

Even if for the function $f$ in \eqref{eq:f} the NHIM remains fixed,
when the local stable and unstable manifolds are extended to global
ones it is expected that in general they will no longer coincide and
indeed they will intersect transversally along a homoclinic manifold. The main tool to study the
splitting of the separatrix \eqref{eq:wsu} as well as the associated
scattering map is the \emph{Melnikov potential} associated to a
perturbation $h$ and to the homoclinic orbit $(p_0,q_0)$:
\begin{equation}\label{eq:melnikov}
\begin{array}{rcl}
\mathcal{L}(I,\varphi,s) & = & \pbe -\int_{-\infty}^{+ \infty} (h(p_0(\sigma),q_0(\sigma),I,\varphi+I \sigma,s+ \sigma;0) \\
&  & \pbe \qquad \quad - h(0,0,I,\varphi+I \sigma,s+\sigma;0))d \sigma, \\
\end{array}
\end{equation}
which taking into account the expression \eqref{eq:perturbation} for
$h$, takes the form
\[
\LL(I,\vp,s)=\int_{-\infty}^{\infty}
[f(q_0(\sigma))-f(0)]g(\vp+I\sigma,s+\sigma) d \sigma.
\]
In our concrete case $f(q)=\cos q$ of \eqref{eq:f}, the Melnikov potential turns out
to be
\[
\LL(I,\vp,s)=\frac{1}{2}\int_{-\infty}^{\infty}
p_{0}^{2}(\sigma)g(\vp+I\sigma,s+\sigma) d \sigma,
\]
and the integral can be explicitly computed by the residue theorem:
\begin{equation}\label{eq:mpex}
\LL(I,\vp,s)=\sum_{(k,l) \in \nn^2 }A_{k,l}(I)\cos(k\vp - ls -
\sigma_{k,l}),
\end{equation}
with \begin{equation}\label{eq:akl}
A_{k,l}(I)=2\pi\frac{(kI-l)}{\sinh\frac{\pi}{2}(kI-l)}a_{k,l},
\end{equation}
where $a_{k,l}$ are the general coefficients of the function $g$
given in \eqref{eq:g}. Notice that the Melnikov potential $\eqref{eq:mpex}$ has
exactly the same harmonics as the perturbation $g$ in \eqref{eq:g}.

We now recall the role played by the Melnikov potential in the
splitting of the separatrix \eqref{eq:wsu}. By Proposition 9.2 in
\cite{DelshamsLS06}, for any $(I,\vp ,s)\in [-1/2,I^{*}_{+}] \times
\tor^{2}$ and for any non-degenerate critical point
$\tau^{*}=\tau^{*}(I,\vp,s)$ of
\begin{equation}\label{eq:tau0}
\tau \in \rr \mapsto \LL (I,\vp-I\tau,s-\tau)
\end{equation}
there exists a locally unique point $z$,
\begin{equation}\label{eq:zhom}
z=z(I,\vp,s;\eps)=(p_0(\tau^{*}),q_0(\tau^{*}),I,\vp,s) +
\ocal(\eps),
\end{equation}
such that $z \in W^s(\tilde \Lambda_{\eps}) \pitchfork W^u(\tilde
\Lambda_{\eps})$.

Next, we are going to find open sets of
$(I,\vp ,s)\in [-1/2,I_{+}^{*}] \times \tor^{2}$, such that the
function \eqref{eq:tau0} has non-degenerate critical points at $\tau
= \tau^{*}(I,\vp,s)$. 


Taking into account the explicit expression for the Melnikov potential \eqref{eq:mpex}, the function \eqref{eq:tau0} takes the form
\begin{equation}\label{eq:Lips}
\LL(I,\vp-I\tau,s-\tau)= \sum_{(k,l) \in \nn^2 } A_{k,l}(I) \cos (k
\vp - ls - \tau(kI-l)),
\end{equation}
with $A_{k,l}(I)$ as in \eqref{eq:akl}. Notice that the Fourier
coefficients $A_{k,l}(I)$ are nothing else but the Fourier
coefficients $a_{k,l}$ multiplied by a non-zero factor depending on
$kI-l$ (which decreases exponentially in $|kI-l|$ as $|kI-l|$ goes to infinity).

The main reason for the introduction of the upper bounds for $|a_{k,l}|$ in \eqref{eq:aklb1} is to make all the computations
for the series defining $\LL(I,\vp,s)$ in \eqref{eq:mpex} and \eqref{eq:Lips} in terms of $\LL^{[\leq 1]}(I,\vp,s)$,
its first order trigonometric polynomial in the angles $(\vp,s)$. Thus, we have
\begin{eqnarray}
\LL(I,\vp,s) & = & A_{0,0} + A_{1,0}(I) \cos \vp + A_{0,1} \cos s + \ocal_2(\rho, r) \nonumber \\
 & := & \LL^{[\leq 1]}(I,\vp,s) + \LL^{[> 1]}(I,\vp,s),\label{eq:Lex}
\end{eqnarray}
where $A_{0,0}=4a_{0,0}$, 
\begin{equation}\label{eq:a01}
A_{0,1}=\frac{2
\pi}{\sinh(\pi/2)}a_{0,1}, \quad \textrm{and} \quad
A_{1,0}(I)=\frac{2 \pi I}{\sinh(\pi/2 I)}a_{1,0}. 
\end{equation}

In the formula above, without lose of genericity and to avoid cumbersome
notation and shifts in the pictures, we have assumed that
$\sigma_{1,0}=\sigma_{0,1}=0$. Otherwise, we can just make a
translation in the variables $(\vp,s)$.

Next we will make our computations for the function
$\LL^{[\leq 1]}$ and a posteriori we will justify that they are also
valid for the complete function $\LL$.

So fixing $(I,\vp,s) \in [-1/2,I_{+}^{*}] \times
\tor^2$ we only need to study the evolution of $\LL^{[\leq 1]}$
along the straight lines
\begin{equation}\label{eq:line}
R: \tau \in \rr \mapsto (\vp-I \tau,s-\tau) \in \tor^2
\end{equation}
on the torus.

\begin{figure}
 \begin{center}
 \includegraphics[height=80mm]{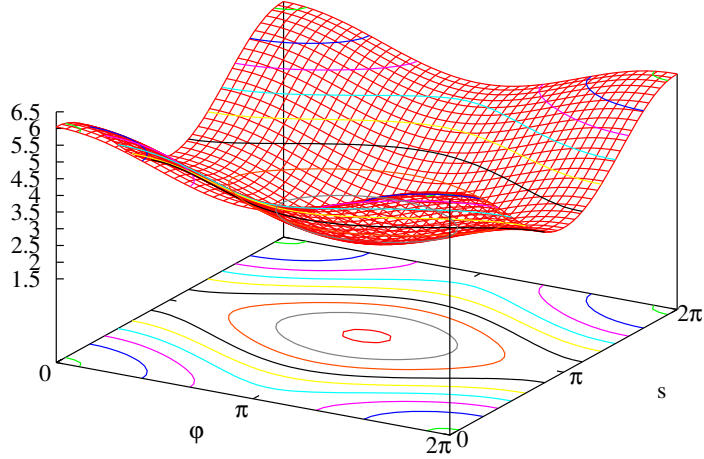}
\caption{Graph and level curves of the Melnikov potential
$\mathcal{L}^{[\leq 1]}(I,\varphi,s)$ with $a_{1,0}=1/4$,
$a_{0,1}=1/2$ and $I=1$. In this case, $A_{0,0}=4$, $A_{1,0}(1)=
\pi/(2\sinh(\pi/2))$ and $A_{0,1}=\pi/\sinh(\pi/2)$}
\label{fig:fig2}.
\end{center}
\end{figure}

By hypothesis \eqref{eq:aklb1}, $a_{1,0} \neq 0$ and $a_{0,1} \neq 0$, and therefore
$|A_{0,1}| \neq 0$ and $|A_{1,0}(I)| \neq 0$ for any $I$. Consequently, for every fixed $I$, the first order trigonometric polynomial 
$(\vp,s) \mapsto \LL^{[\leq 1]}(I,\vp,s)$ possesses four
non-degenerate critical points at $(0,0)$, $(0,\pi)$, $(\pi,0)$ and
$(\pi,\pi)$; a maximum, a minimum and two saddle points. Without loss of generality and for illustration purposes let us assume $a_{1,0}, a_{0,1}
>0$, so that $A_{1,0}(I), A_{0,1}>0$ for any $I$. In this way the maximum of $\LL^{[\leq 1]}(I,\cdot,\cdot)$ is attained at $(0,0)$, the minimum at $(\pi,\pi)$ and the two saddles at
$(0,\pi)$ and $(\pi,0)$ (see Figure \ref{fig:fig2}). Of course, assuming that
$0<\rho \leq \rho^{*}$ and $0<r \leq r^{*}$, for $\rho^{*}$ and $r^{*}$ small enough, by the implicit function Theorem, 
the function $\LL(I,\cdot,\cdot)$ possesses also four non-degenerate critical points close to ones of  $\LL^{[\leq 1]}(I,\cdot,\cdot)$, with the same properties.

Around the two extremum points (the maximum and the minimum), its level curves are closed curves which fill out a basin ending at the level curve of one of
the saddle points. Therefore, any straight line \eqref{eq:line} that enters into one of the two 
extremum basins is tangent to one of the closed level curves,
giving rise to one or more extrema of \eqref{eq:tau0}.
Since the two extrema of $\LL(I,\cdot,\cdot)$ are non-degenerate for any $I$, the closed level curves close to them are convex. Therefore,
any straight line \eqref{eq:line} passing close enough to the extrema gives rise to unique non-degenerate extremum of \eqref{eq:tau0}. In particular, for irrational values of
$I$ the line \eqref{eq:line} is a dense straight line, so there exist an infinite number of non-degenerate extrema for \eqref{eq:tau0}. Nevertheless, and independently
of the irrational character of $I$, thanks to the form
of the perturbation, we are going to see that indeed \emph{all} the closed level curves in any extremum basin are convex.


To do so, we proceed in the following way.  Given a fixed value of
$I$, let us look for the geometric locus where the straight lines \eqref{eq:line} are tangent to the level curves of $\LL(I,\cdot,\cdot)$.
For this, we have to impose that $\nabla_{\vp,s} \LL$, the gradient of $\LL(I,\cdot,\cdot)$,
is orthogonal to the slope $(I,1)$ of the straight line \eqref{eq:line}:
\begin{equation}\label{eq:implicit}
I \frac{\partial{\LL}}{\partial \vp} (I,\vp,s) +
\frac{\partial{\LL}}{\partial s} (I,\vp,s) =0.
\end{equation}

Intuitively, for fixed $I$, if we want to pass through a mountain of height $\LL(I,\vp,s)$ along straight lines following a direction $(I,1)$,
equation \eqref{eq:implicit} gives the position of
the points  $(\vp,s)$ of maximum height, the \emph{crest}, that we denote by $\mathcal{C}=\mathcal{C}(I)$.

Using the expression \eqref{eq:Lex} for $\LL$, equation \eqref{eq:implicit} has the form
\[ I A_{1,0}(I) \sin \vp + A_{0,1} \sin s + \ocal_2(\rho,r) = 0.\]

Disregarding first the $\ocal_2(\rho,r)$ term we are faced with the following implicit equation
\begin{equation}\label{eq:implicits}
\alpha(I) \sin \vp + \sin s=0, 
\end{equation}
where
\begin{equation}\label{eq:alpha}
\alpha(I):= \frac{I A_{1,0} (I)}{A_{0,1}} = \frac{\sinh(\pi/2)
I^2}{\sinh(\pi/2 I)}\frac{a_{1,0}}{a_{0,1}}.
\end{equation}

\begin{figure}
 \begin{center}
 \includegraphics[height=100mm]{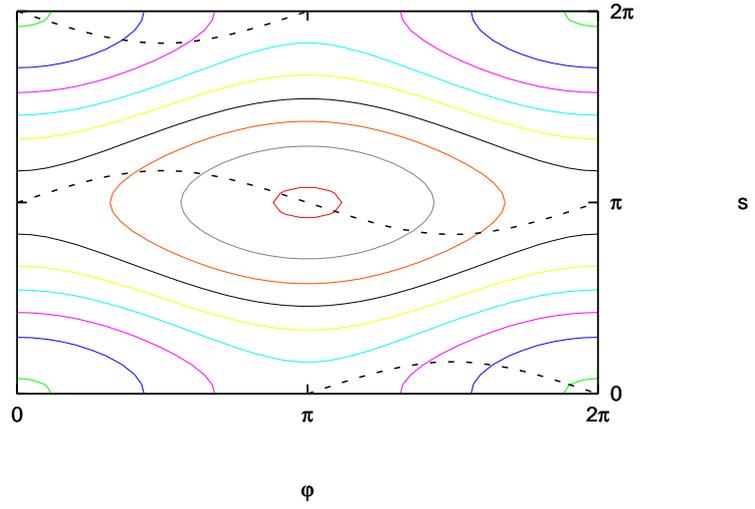}
 \caption{Closed curves satisfying \eqref{eq:implicits} for $I=1$ (the crests), dashed
black, and level sets of the function $\LL^{[\leq 1]}(1,\cdot,\cdot)$ with the same values as
in Figure \ref{fig:fig2}.}
 \label{fig:implicit}
 \end{center}
\end{figure}

Assuming that
\[ |\alpha(I)| <1, \]
which holds for all $I$ as long as
\begin{equation}\label{eq:b}
1.03 \left | \frac{a_{1,0}}{a_{0,1}}\right | < 1,
\end{equation}
and that $\rho$ and $r$ are small enough, equation \eqref{eq:implicit} defines exactly two closed curves $\mathcal{C}_M$ and
$\mathcal{C}_m$, parameterized by $\vp$, which will be called \emph{crests}. The crest
$\mathcal{C}_M=\mathcal{C}_M(I)$, passing through the maximum
$(0,0)$ of $\LL$ contains the saddle $(0,\pi)$ and is given by the expression
\begin{equation}\label{eq:curvemax}
s=\xi_M(\vp,I)=- \arcsin ( \alpha(I) \sin \vp) + \ocal_2(\rho,r).
\end{equation}
The crest $\mathcal{C}_m=\mathcal{C}_m(I)$ passes through the
minimum $(\pi,\pi)$ and contains the saddle $(\pi,0)$. It is given by the expression
\begin{equation}\label{eq:curvemin}
s=\xi_m(\vp,I)=-\arcsin ( \alpha(I) \sin (\vp + \pi)) + \pi + \ocal_2(\rho,r).
\end{equation}
In Figure \ref{fig:implicit} there appear these two curves (dashed
black) as well as the level sets of the function $\LL^{[\leq 1]}$.

\begin{rem}
The case when $|\alpha(I)|>1$ is analogous, but then the crests
$\mathcal{C}_m$ and $\mathcal{C}_M$ are parameterized by the
variable $s$. The case $|\alpha(I)|=1$ is special, because in this
case the union of the two curves degenerates in two straight lines
along which the function $\LL$ is constant.

These two additional cases do not present additional difficulties for the
construction of the function $\tau^{*}$ and therefore for the
subsequent mechanism of diffusion. However, for easiness of the
reading we are going to concentrate simply in the case \eqref{eq:b},
which is fulfilled thanks to the first inequality of hypothesis \eqref{eq:clambda},
which is the only one that we are going to assume from now on.
\end{rem}

For any point $(I,\vp,s) \in [-1/2,I_{+}^{*}] \times \tor^2 \times \mathbb{R}$, the value $\tau^{*}$ for which
the function \eqref{eq:tau0} has a non-degenerate critical point
satisfies $(I,\vp-I \tau^{*}, s- \tau^{*}) \in \mathcal{C}_m \cup
\mathcal{C}_M$. Thus, the non-degenerate critical points $\tau^{*}$ of \eqref{eq:tau0} are achieved at the intersection of the straight line \eqref{eq:line} with either
the crest $\mathcal{C}_M$ or $\mathcal{C}_m$.

\begin{figure}
 \begin{center}
 \begin{tabular}{cc}
 \includegraphics[width=60mm]{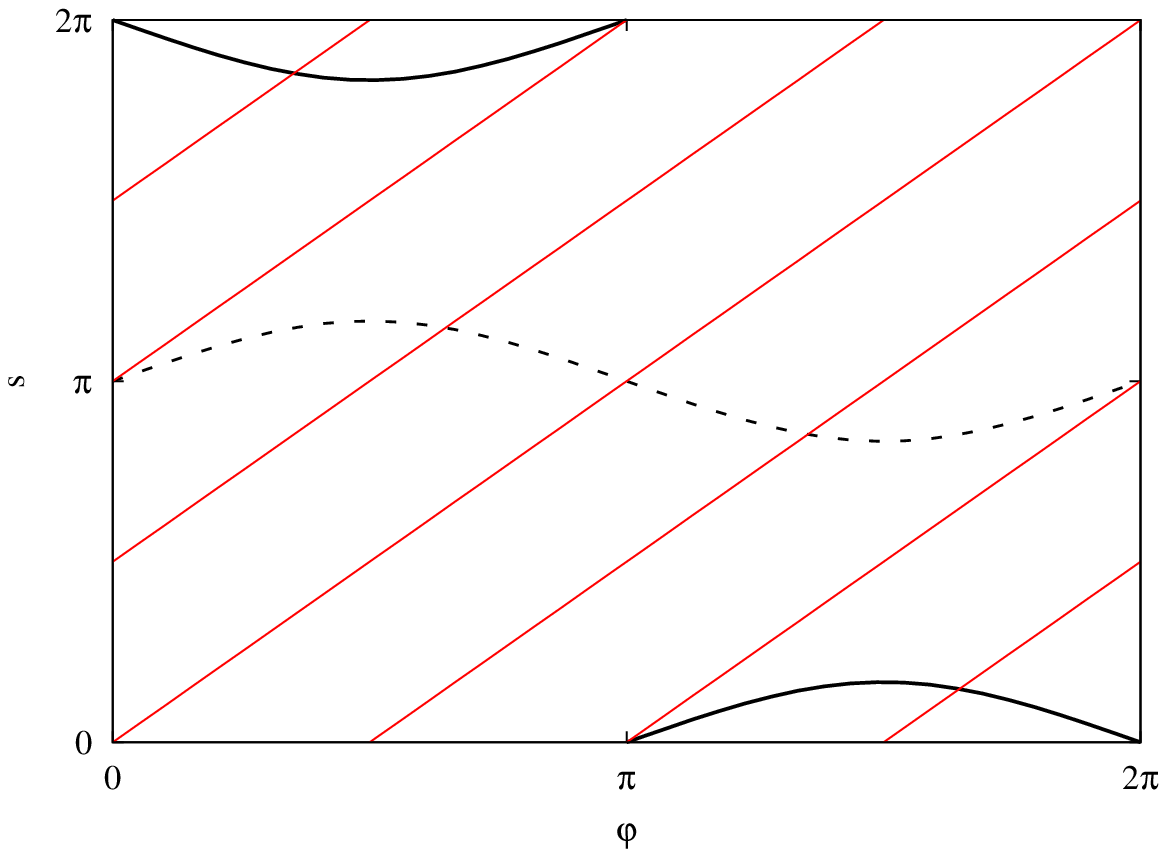} &
 \includegraphics[width=60mm]{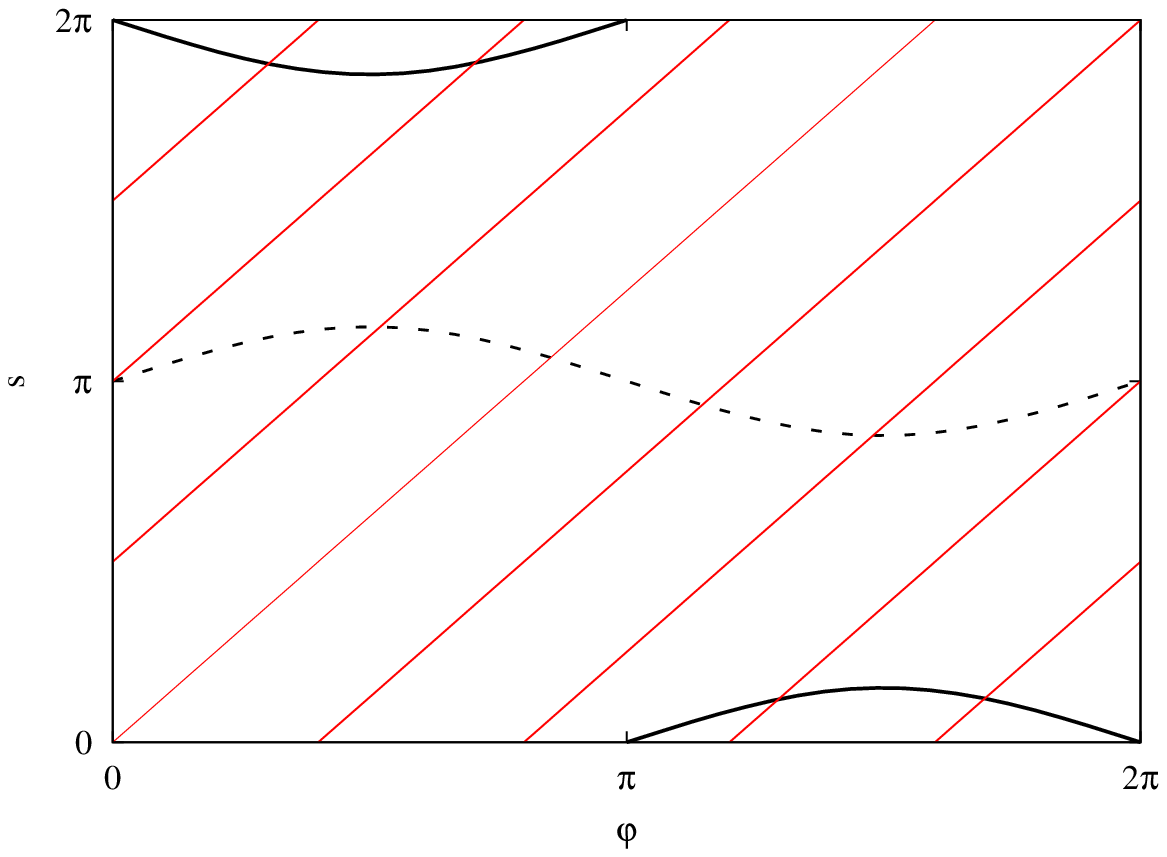}\\
 \end{tabular}
 \caption{Straight lines \eqref{eq:line} with slope 1 (Left) and 0.8 (Right) and the curves
 \eqref{eq:implicits}:  the curve of the maxima (solid curve) and the curve of the minima (dashed curve). See the text. }
 \label{fig:fig3}
 \end{center}
\end{figure}

Of course, for any point $(I,\vp,s)$, there may exist several
intersections of the line \eqref{eq:line} with the crests
$\mathcal{C}_M$ and $\mathcal{C}_m$, giving rise to different
homoclinic intersections and (as we will see in the next section)
different scattering maps. See Figure \ref{fig:fig3}.

From now on, we choose only one of these intersections, the one with the crest $\mathcal{C}_M$ and we define the
function $\tau^{*}$ as the one that given $(I,\vp,s)$,
associates the minimum $|\tau|$ such that
\[ (I,\vp-I\tau,s-\tau) \in \mathcal{C}_M(I).\]
In symbols,
\begin{equation}\label{eq:taudef}
\tau^{*}(I,\vp,s):= \min_{|\tau| \in \mathbb{R}} \{ (I,\vp-I\tau,s-\tau) \in \mathcal{C}_M(I) \}. 
\end{equation}

To determine a domain of definition as large as possible where
this function is continuous and to avoid a casuistic discussion, we need to check that for any fixed $I$, the straight lines cross only
once the crest $\mathcal{C}_M$ inside the domain $(\vp,s)$ and they
do it transversally. This implies that for any fixed $I$, the slope
$1/I$ of the straight lines is bigger than the slope of the
derivative of the function $\xi_M(\vp,I)$ for all $\vp \in \tor$,
that is
\begin{equation}\label{eq:transv}
\frac{1}{I} > \frac{\partial \xi_M}{\partial \vp}(\vp,I), \quad \mathrm{for} \, \,  \mathrm{all} \, \vp \in \tor,
\end{equation}
which by equations
\eqref{eq:curvemax} and \eqref{eq:alpha} and assuming that $\rho$ and $r$ are small enough, is equivalent to
\[ \alpha (I) I <1.\]
By expression \eqref{eq:alpha} for $\alpha(I)$, it is easy to see that
\[ \alpha (I) I < 1.6 \frac{a_{1,0}}{a_{0,1}},\]
so that condition \eqref{eq:transv} holds for all $I$ as long as
\begin{equation}\label{eq:coefcond}
1.6 \left | \frac{a_{1,0}}{a_{0,1}} \right | < 1,
\end{equation}
which is exactly hypotheses \eqref{eq:clambda} of Theorem \ref{teo:main} and implies \eqref{eq:b}.

Under condition \eqref{eq:coefcond}, one possible domain of
definition $H$ where $\tau^{*}$ is continuous consists of excluding, for
any $I \in [-1/2,I^{*}_{+}]$ the crest $\mathcal{C}_m(I)$ from the domain
of $(\vp,s)$, that is
\begin{equation}\label{eq:H2dom}
H=\{ (I,\vp,s) \in [-1/2,I_{+}^{*}] \times \tor \times \mathbb{R}: \xi_m(I,\vp) -2 \pi < s < \xi_m(I,\vp)\}.
\end{equation}
See Figure \ref{fig:domain}.

This particular choice of $\tau^{*}$ gives rise to a homoclinic manifold to which we will associate in the next
section an outer dynamics to the NHIM, that we will describe by means of the scattering map.

\begin{rem}
Condition \eqref{eq:coefcond} is very convenient since it provides a large domain of definition $H$ for $\tau^{*}$, and therefore it allows us to define a global
homoclinic manifold $\Gamma_{\eps}$. Although the condition is very restrictive, it is 
not necessary for the mechanism of diffusion. Indeed, we introduced this assumption just to avoid a casuistic description. However, if the condition is not satisfied
we can obtain several homoclinic manifolds giving rise to different scattering maps, offering more possibilities for diffusion. It also opens the field for 
homoclinic bifurcations for a NHIM (see \cite{DelshamsMR08}). 
\end{rem}

\begin{rem}\label{rem:Lalongc}
The crest $\mathcal{C}_M$ is in the maximum basin of $\LL$, so the
function $\LL$ decreases when one travels from $(0,0)$ to
$(0,\pi)$ increasing $\vp$ along the curve $\mathcal{C}_M$ and
increases when one travels from $(0,\pi)$ to $(0,2 \pi)$.
Analogously, the other crest $\mathcal{C}_m$ is in the minimum
basin of $\LL$ with a complete analogous property. Since by \eqref{eq:implicit} and
\eqref{eq:transv}, the curve $(\vp,\xi_M(I,\vp))$ is never tangent
to the level sets of the function $\LL$, the increase and decrease are \emph{strict}.
\end{rem}


\subsection{Part 2. Outer dynamics (Scattering
map)}\label{sec:part2}

We have seen in the previous section that if condition \eqref{eq:clambda} is satisfied
for any $(I,\vp,s)$ in the domain $H$ given in \eqref{eq:H2dom}, the function \eqref{eq:tau0} has a
non-degenerate critical point $\tau^{*}$ given by
$\tau^{*}=\tau^{*}(I,\vp,s)$, where $\tau^{*}$ is a smooth function
defined in \eqref{eq:taudef}. If $0<\eps < \eps^{*}(I^{*}_{+})$, these
critical points give rise to a homoclinic manifold
$\Gamma_{\eps} \subset W^{s}_{\tilde \Lambda}  \pitchfork W^{u}_{\tilde \Lambda}$,
along which these invariant manifolds intersect
transversally. By equation \eqref{eq:zhom}, it has the form
\begin{equation*}
\begin{array}{r@{}l}
\Gamma_{\eps} = \{ & z =
z(I,\vp,s;\eps)=(p_0(\tau^{*}),q_0(\tau^{*}),I,\vp,s) + \ocal(\eps): \\ 
& (I,\vp,s) \in H, \tau^{*}=\tau^{*}(I,\vp,s) \in \mathbb{R}\}.
\end{array}
\end{equation*}

\begin{rem}
For the experts in the splitting of separatrices, we notice that the size of $\eps^{*}$ required for the justification of the transversal intersection of 
$W^{s}_{\tilde \Lambda}$ and $W^{u}_{\tilde \Lambda}$ along $\Gamma_{\eps}$ has to be such that the Melnikov potential \eqref{eq:mpex} gives the dominant part of the 
formula for the splitting. In our case, since $\LL$ as well as its two first derivatives are $\ocal\left(\exp\left(-\pi/2 I^{*}_{+}\right)\right)$
on the domain $H$, we need to impose that $\eps^{*}=\ocal\left(\exp\left(-\pi/2 I^{*}_{+}\right)\right)$.
\end{rem}

In this section we will see that associated to the homoclinic manifold
$\Gamma_{\eps}$ we can define an outer dynamics $S_{\eps}$ to the NHIM
$\widetilde \Lambda$ and we will obtain an approximate explicit
expression for it.

The scattering map associated to $\Gamma_{\eps}$ is defined in the
following way:
\begin{equation}\label{eq:scat}
\begin{array}{rlcc}
S_{\eps}: & H \subset \widetilde{\Lambda} &  \longrightarrow & \widetilde{\Lambda} \\
 & x_{-} & \mapsto & x_{+} \\
\end{array}
\end{equation}
such that $x_{+}=S_{\eps}(x_{-})$ if and only if there exists $z \in
\Gamma_{\eps}$ such that
\[ \mathrm{dist}(\Phi_{t,\eps}(z),\Phi_{t, \eps}(x_{\pm})) \rightarrow 0 \qquad \mathrm{ for } \, t \rightarrow \pm \infty,\]
where $\Phi_{t,\eps}$ is the flow of Hamiltonian \eqref{eq:example}.

In words, the scattering map maps a point $x_{-}$ on the NHIM to a point $x_{+}$
on the NHIM if there exists a homoclinic orbit to the NHIM
that approaches the orbit of $x_{-}$ in the past and the orbit of
$x_{+}$ in the future.


The scattering map $S_{\eps}$ is exact and symplectic and indeed it is Hamiltonian, since it is given by the
time $\eps$ map of a Hamiltonian $\mathcal{S}_{\eps}$ \cite{DelshamsLS08}. In the variables $(I,\vp,s)$ this implies that the 
following formula holds for the scattering map
\begin{equation*}
S_{\eps}(I,\vp,s)= (I,\vp,s) + \eps \left ( - \frac{\partial
\mathcal{S}_0}{\partial \vp} (I,\vp,s) , \frac{\partial
\mathcal{S}_0}{\partial I} (I,\vp,s),0 \right ) + \ocal (\eps^2).
\end{equation*}

As it is described in equation (21) in \cite{DelshamsH09}, the dominant term $\mathcal{S}_0$ of the Hamiltonian $\mathcal{S}_{\eps}$ is equal with opposite sign to the 
\emph{reduced Poincar\'e function} $\LL^{*}$ defined implicitly by
\begin{equation}\label{eq:rpf}
\LL^{*}(I,\vp - Is):=\LL(I,\vp - I \tau^{*}(I,\vp,s),
s-\tau^{*}(I,\vp,s)).
\end{equation}

It is important to notice that the Hamiltonian $\mathcal{S}_0$ is an autonomous $1$-degree of freedom Hamiltonian that can be expressed in terms of the variable 
$\tet = \vp - I s$:
\begin{equation}\label{eq:hamscat}
\mathcal{S}_{0}(I,\vp,s)=-\LL^{*}(I,\tet), \quad  \tet=\vp-Is,
\end{equation}
so that
\begin{equation}\label{eq:scatt}
S_{\eps}(I,\vp,s)= \left ( I + \eps \frac{\partial
\LL^{*}}{\partial \tet} (I,\tet) + \ocal (\eps^2), \tet - \eps \frac{\partial \LL^{*}}{\partial
I} (I,\tet) + \ocal (\eps^2), s \right ),
\end{equation}
and the iterates under the scattering map simply follow closely the level curves of the reduced Poincar\'e function \eqref{eq:rpf}.

In order to obtain an expression for the reduced Poincar\'e function and for its level curves in our particular example, we will perform a discussion based on
geometric considerations.

By the definition of $\tau^{*}(I,\vp,s)$ given in the previous
section we have that the point
\begin{equation}\label{eq:cpoint}
c(I,\vp,s):=(I,\vp - I \tau^{*}(I,\vp,s),s-\tau^{*}(I,\vp,s)) \in \mathcal{C}_M(I),
\end{equation}
belongs to the crest $\mathcal{C}_M$, which is the closed curve defined in \eqref{eq:curvemax}.

\begin{figure}
 \begin{center}
 \includegraphics[width=80mm]{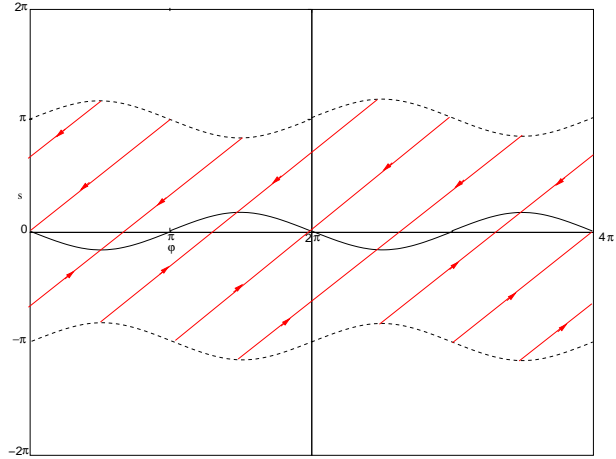}
 \caption{Straight lines \eqref{eq:line} with slope $1$ and the curves
 \eqref{eq:implicits}:  the curve of the maxima (solid curve) and the curve of the minima (dashed curve). See the text. }
 \label{fig:domain}
 \end{center}
\end{figure}

Therefore, the reduced Poincar\'e function $\LL^{*}$ evaluated on a
point $(I,\vp,s)$ in the domain $H$ defined in \eqref{eq:H2dom} provides the value of the function $\LL$
evaluated on $c(I,\vp,s)$, the closest intersection of the straight line
\eqref{eq:line} starting on this point $(I,\vp,s)$ with the curve
$\mathcal{C}_M$. By construction, it is clear that there is a segment of points $(I,\vp,s)$ in the domain $H$ with the same $c(I,\vp,s)$ on the curve $\mathcal{C}_M$.
See Figure \ref{fig:domain}.


Since the function $\LL^{*}$ is constant on these segments it can be written as a function of only two variables: the action $I$ and the variable $\tet=\vp - I s$, which is
$2\pi$-periodic in $\vp$ and constant along the straight lines \eqref{eq:line} of slope $1/I$
contained in \eqref{eq:H2dom}.


In order to obtain explicitly the expression of $\LL^{*}$ in the
variable $\tet$, we will proceed in the following way. Fixed $I$, we consider a point $(I,\vp^{*},0)$ lying on the axis $s=0$. The points $(I,\vp,s)$ in the domain $H$ 
defined in \eqref{eq:H2dom} satisfying $\tet=\vp-Is=\vp{*}$ intersect the 
curve $\mathcal{C}_M$ on the same point $c(I,\vp^{*},0)$ as $(I,\vp^{*},0)$ which by \eqref{eq:cpoint} is 
\[(I,\vp^{*}-I\tau^{*}(I,\vp^{*},0),-\tau^{*}(I,\vp^{*},0)).\] 
Since each segment $\{\tet=\cte \}$ in $H$ intersects the $s=0$ axis at the point $(I,\tet,0)$, 
the function $\LL^{*}$ has the following expression:
\[ \LL^{*}(I,\tet)=\LL(I,\tet-I
\tau^{*}(I,\tet,0),-\tau^{*}(I,\tet,0)).\]

The behavior of the function $\LL^{*}$ with respect to the variable
$\tet$ (which parameterizes the curve $\mathcal{C}_M$) is exactly the
behavior of the function $\LL$ along the curve $\mathcal{C}_M$, that we discussed in Remark
\ref{rem:Lalongc}. Namely, when $\tet$ increases from $0$ to $\pi$, one travels from $(0,0)$ to $(0,\pi)$, increasing $\vp$ along the curve $\mathcal{C}_m$ and therefore
the function $\LL^{*}$ decreases strictly. Equivalently, when $\tet$ increases from $\pi$ to $2\pi$, one travels from $(0,\pi)$ to $(0,2 \pi)$, increasing $\vp$ along the 
curve $\mathcal{C}_m$ and therefore the function $\LL^{*}$ increases strictly.

We summarize some of the properties of the reduced Poincar\'e
function in the following Proposition:
\begin{prop}\label{prop:prop}
For any $I \in [-1/2,I^{*}_{+}]$, the function $\tet \mapsto
\LL^{*}(I,\tet)$ has a non-degenerate maximum (minimum) close to $\tet=0
\, (\mathrm{mod} 2\pi)$ ($\tet=\pi \, (\mathrm{mod} 2\pi)$, respectively)
and is strictly monotone in all the other points $\tet$. Moreover, it has
the following expression
\begin{eqnarray}
\LL^{*}(I,\tet) & = &  A_{0,0} + A_{1,0} (I) \cos (\tet - I
\tau^{*}(I,\tet,0))
+  A_{0,1} \cos ( \xi_M (I, \tet - I \tau^{*}(I,\tet,0))) \nonumber \\
& & + \ocal_2 (\rho,r), \label{eq:Ls}
\end{eqnarray}
where $A_{1,0}$ and $A_{0,1}(I)$ are given in \eqref{eq:a01}, $\xi_M$ in \eqref{eq:curvemax} and $\tau^{*}(I,\tet,0)$ is defined in \eqref{eq:taudef}.
\end{prop}

\begin{rem}\label{rem:rpf0}
Notice that the behavior of the function $\LL^{*}$ with respect to the variable $\tet$ is ``cosinus-like''. This observation is clear  when one considers the case 
$I=0$, where $\vp=\tet$, $\xi_M(0,\vp)=0$ and
\[
\LL^{*}(0,\tet)  =   A_{0,0} +  A_{1,0} (0) \cos (\tet) + A_{0,1} +
\ocal_2 (\rho,r). \]
\end{rem}

\begin{figure}[ht]
 \begin{center}
 \includegraphics[width=120mm]{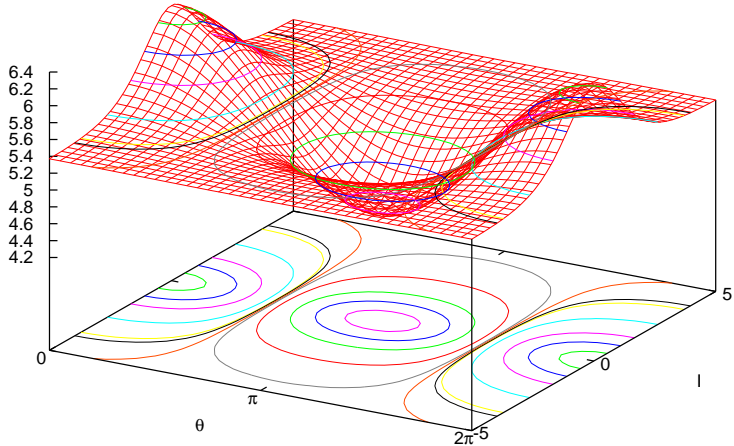}
 \caption{Graph and level curves of the reduced Poincar\'{e} function $\LL^{*}(I,\tet)$ given in \eqref{eq:Ls}. As in
 Figure \ref{fig:fig2} we have chosen $a_{1,0}=1/4$ and $a_{0,1}=1/2$ in \eqref{eq:a01}.}
 \label{fig:levelsets}
 \end{center}
\end{figure}

Proposition \ref{prop:prop} provides us with an exhaustive description of the level sets of the reduced Poincar\'e function, giving an approximate expression in first
order for the orbits of the scattering map $S_{\eps}$ in \eqref{eq:scatt}. In Figure \ref{fig:levelsets} we illustrate these level curves for a particular case.

\subsection{Part 3. Inner dynamics}\label{sec:inner}

The inner dynamics is the dynamics of the flow of Hamiltonian
\eqref{eq:example} restricted to the NHIM \eqref{eq:nhim}. Indeed, by
the form of the perturbation, the Hamiltonian restricted to NHIM
\eqref{eq:nhim} takes the explicit form
\begin{equation}\label{eq:hamin}
K(I,\vp,s;\eps)=\frac{I^2}{2}+ \eps g(\vp,s),
\end{equation}
where $(I,\vp,s) \in [-1/2, I_{+}^{*}] \times \tor^2$ and $g$ is given in
\eqref{eq:g}.

In this section, we want to study the dynamics in the NHIM, that is,
what are the invariant objects, what is the distance among them in
terms of the action $I$ and what are their approximate analytical
expression. This section relies on the proof and the results of
Theorem 3.1 in \cite{DelshamsH09}.

We already mentioned in Section \ref{sec:nhim} that the dynamics for
the unperturbed case (that is, $\eps=0$) is very simple. Indeed, all
the trajectories lie on two-dimensional invariant tori
$I=\mathrm{const}$. The motion on the tori is conjugate to a rigid
rotation of frequency vector $(I,1)$. Notice that the Hamiltonian
$K$ is non degenerate, that is,
\[\frac{\partial^2 K}{\partial I^2}
\equiv 1 \neq 0.\]

For $\eps>0$, KAM theorem ensures the preservation, with some
deformation, of most of the invariant tori present in the
unperturbed system. Indeed, those having frequencies 
``sufficiently'' non-resonant, for which the influence of small
divisors $kI-l$ can be overcome. This is guaranteed by a Diophantine
condition on the frequency vector:
\[ |kI-l| \geq \frac{\gamma}{|(k,l)|^{\tau}} \quad \forall (k,l) \in \mathbb{Z}^2 \backslash \{0\}, \]
with $\tau \geq 1$ and some $\gamma=\ocal(\sqrt{\eps})>0$. The
frequency vectors $(I,1)$ satisfying this Diophantine condition fill
a Cantorian set of relative measure $1-\ocal (\sqrt{\eps})$, called the \emph{non-resonant region}.

Hence, the invariant tori with Diophantine frequencies persist under
the perturbation, with some deformation. These KAM tori, which are
just a continuation of the ones that existed in the integrable case
$\eps=0$, are commonly known as \emph{primary tori} and are
given by the level sets of a function $F$ of the form (see
Proposition 3.24 in \cite{DelshamsH09})
\[F(I,\vp,s) = I + \ocal(\eps).\]

On the contrary, the invariant tori with resonant frequencies are typically
destroyed by the perturbation, creating gaps in the foliation of
invariant tori of size up to $\mathcal{O}(\sqrt{\eps})$ centered
around resonances (indeed for the values of $I$ such that $kI-l=0$,
for some $(k,l) \in \mathbb{N}^2$ which is the support of the
Fourier transform of the perturbation $g$ given in \eqref{eq:g}).
However, in these resonant regions, other invariant objects are
created, like \emph{secondary tori}, that is, two-dimensional invariant KAM tori
contractible to a periodic orbit, as well as periodic orbits with associated stable and unstable manifolds. To
prove the existence of these objects in the \emph{resonant regions}
and also to give an approximate expression for them, in
\cite{DelshamsLS06, DelshamsH09} several steps of averaging were
performed before applying the KAM theorem to the Hamiltonian
expressed in the averaged variables.

More precisely, given any $(k_0,l_0) \in \nn^2$, $k_0 \neq 0$,
$\gcd(k_0,l_0) =1$, for a resonant region centered around a resonance
$I=l_0/k_0$, the invariant tori are given by the level sets
of a function $F$, whose dominant term $\bar F$ in $\eps$ is of the
form (see Theorem 3.28 in \cite{DelshamsH09})
\begin{equation}\label{eq:Forder}
\bar F (I,\theta) = \frac{(k_0 I - l_0)^2}{2} + \eps k_0^2
U^{k_0,l_0}(\theta),
\end{equation}
where $\theta=k_0 \vp - l_0 s$ and the function $U^{k_0,l_0}$
contains the resonant terms of $g$
with respect to $(k_0,l_0)$, that is,
\begin{equation}\label{eq:Uex}
U^{k_0,l_0}(\theta)  =  \sum_{t=1}^{M}  a_{tk_0,tl_0} \cos (t
\theta),
\end{equation}
where $M=\eps^{-1/(26+\delta)}$ and $\delta$ is a parameter
satisfying  $0<\delta<1/10$. The lower and upper bounds for the coefficients $a_{kl}$
of $g$ provided by hypothesis \eqref{eq:aklb1} ensure that the function
\eqref{eq:Uex} reaches a non-degenerate global maximum for any
$I=l_0/k_0$ and any $M \geq 1$, independently of $\eps$. 

Indeed, the function $U^{k_0,l_0}(\theta)$ is analytic and $2
\pi$-periodic in $\theta$. By \eqref{eq:aklb1}, the function
$U^{k_0,l_0}$ in \eqref{eq:Uex} is well approximated by its first order trigonometric polynomial
\begin{eqnarray}\label{eq:Uaproxb}
U^{k_0,l_0}(\theta) & = &  a_{k_0,l_0} \cos (\theta)  + \ocal_2
(\rho^{k_0} r^{l_0}).
\end{eqnarray}
As long as $\rho,r$ are small enough, and replacing $\theta$ by $\theta+\pi$ if necessary (when
$a_{k_0,l_0}<0$), 
the function \eqref{eq:Uaproxb} has two non-degenerate critical points corresponding to a global maximum
at $\theta=0$ and a global minimum at $\theta=\pi$.

Notice that in the resonant region around $I=l_0/k_0$, the angle variable $\tet$ introduced in
\eqref{eq:hamscat} satisfies $\tet=\vp-(l_0/k_0) s =\theta/k_0$.
Hence, the invariant tori in the resonant region can be approximated
by the level sets of a function that expressed in the same variables
$(I,\tet)$ as the scattering map in \eqref{eq:scatt}, takes the form
\begin{equation}\label{eq:Fsim}
F ^{*}(I,\tet) = \frac{\left ( I - \frac{l_0}{k_0} \right
)^2}{2} + \eps a_{k_0,l_0} ( \cos (k_0 \tet) - 1),
\end{equation}
which is the Hamiltonian of a pendulum in the variables
$(I,\theta=k_0 \tet)$. Notice that a constant term $-\eps
a_{k_0,l_0}$ has been added in order that the 0-level set
corresponds to the separatrices of the pendulum. It is worth
noticing that the size in the action $I$ of the region enclosed by
the two separatrices of \eqref{eq:Fsim}, the gap, is given
by $\sqrt{\eps |a_{k_0,l_0}|}$. In terms of the variable $\tet$, the function $F^{*}$ is $2 \pi k_0$-periodic and therefore 
the region enclosed by the separatrices has $k_0$ components, the
``eyes''.

The result of Theorem 3.1 in \cite{DelshamsH09} provides a sequence
of KAM tori consisting of primary and secondary tori which are
$\eps^{1+\eta}$-close spaced in terms of the action variable $I$,
for some $\eta>0$.

The scattering map \eqref{eq:scat} detects the existence of
heteroclinic connections among primary or secondary invariant tori
whose distance is smaller than $\eps$. Thus, we will distinguish
two types of resonant regions depending whether the size of the gaps
created by the resonances is bigger or smaller than the size $\eps$
of the heteroclinic jumps. The resonant region with big gaps
corresponds to those resonances centered around $I=l_0/k_0$ such
that $\sqrt{|a_{k_0,l_0}|} \geq \eps^{1/2}$ and for this region we
have the \emph{large gap problem}. In these regions, we will include
in the transition chain secondary tori given approximately by
negative level sets of the function \eqref{eq:Fsim}. In the regions
with small gaps, centered around $I=l_0/k_0$ such that
$\sqrt{|a_{k_0,l_0}|} < \eps^{1/2}$, the gaps are of size smaller than
$\eps$ in terms of the $I$ variable, so that it is possible to
connect two primary tori on both sides of the gap. This case does
not present the large gap problem and can be treated analogously as in
the non-resonant region. Both regions, the small gaps region and the
non-resonant region, are called in \cite{DelshamsH09} \emph{the flat
tori region} and the dominant term in $\eps$ of the invariant tori
is given there in first order by the function
\begin{equation}\label{eq:Fordern}
F^{*}(I,\tet)=I,
\end{equation}
where $\tet=\vp - I s$.

%

\subsection{Part 4. Combination of both dynamics}

The geometric mechanism of diffusion close to the NHIM is based on
the combination of two types of dynamics, the inner one, provided by
Hamiltonian \eqref{eq:hamin} and the outer one, approximately given by the time-$\eps$ map of the Hamiltonian \eqref{eq:hamscat}.

Diffusion inside the three-dimensional NHIM can only take place 
between two-dimensional invariant tori. In order to overcome the
obstacles of these invariant KAM tori present in the NHIM, we need to
use the outer dynamics to ``jump'' from one KAM invariant torus to another one. Of course, for this mechanism to be
successful the outer dynamics does not have to
preserve the invariant KAM tori of the inner one.

We have seen that the invariant KAM tori for  Hamiltonian
\eqref{eq:hamin} are given by the level sets of a function $F$,
having different expressions in the flat tori region and in the big gaps region. Indeed, its dominant term is given by a function $F^{*}$ that
depends only on the variables $(I,\tet)$, having expression  \eqref{eq:Fsim} for the big gaps region and expression
\eqref{eq:Fordern} for the flat tori region. Moreover, in equation
\eqref{eq:scatt} we showed that the scattering map $S_{\eps}$ is given in
first order by the time-$\eps$ map of an integrable Hamiltonian
$-\LL^{*}(I,\tet)$ of the form \eqref{eq:Ls}. Therefore, imposing that
$F^{*}$ is not a first integral
of the Hamiltonian $\LL^{*}$, or equivalently, that the functions $F^{*}$ and $\LL^{*}$ are
functionally independent as variables of $(I,\tet)$:
\begin{equation}\label{eq:pb}
\frac{\partial F^{*}}{\partial I} \dot{I} + \frac{\partial
F^{*}}{\partial \tet} \dot{\tet} = - \frac{\partial F^{*}}{\partial
I} \frac{\partial \LL^{*}}{\partial \tet} + \frac{\partial
F^{*}}{\partial \tet} \frac{\partial \LL^{*}}{\partial I} = \{
F^{*}, \LL^{*}\}  \not \equiv 0,
\end{equation}
guarantees that the KAM tori in the NHIM are not invariant for the
scattering map.

Moreover, it is easy to see that the Poisson bracket $\{F,\LL^{*}\}$ provides,
in first order, the deformation under the scattering map of an
invariant torus given by $F(I,\vp,s)=\cte$, since
\begin{equation}\label{eq:Fscat}
F \circ S_{\eps} = F - \eps \{F,\LL^{*}\} + \ocal (\eps^{2}).
\end{equation}
Considering the torus
\[ \{F(I,\vp,s)=E\},\]
the image under the scattering map satisfies $\{F \circ S^{-1}_{\eps}(I,\vp,s) =E\}$, where $F \circ S^{-1}$ has the following expression:
\begin{equation}\label{eq:Escat}
F \circ S^{-1}_{\eps} = F + \eps \{F,\LL^{*}\} + \ocal (\eps^{2}).
\end{equation}

Using that $F$ can be approximated by its dominant term $F^{*}$, it
is enough to check
\begin{equation}\label{eq:pois}
\{F^{*},\LL^{*}\} \not \equiv 0
\end{equation}
to guarantee that the scattering map moves the invariant tori of the inner
dynamics.

In the flat tori region, $F^{*}$ is independent of $\eps$ and is given in \eqref{eq:Fordern}
so that 
\begin{equation}\label{eq:Lflat}
\{F^{*},\LL^{*}\}(I, \tet) = -
\frac{\partial \LL^{*}}{\partial \tet} (I, \tet),
\end{equation}
whereas in the big gaps region $F^{*}$ has the form \eqref{eq:Fsim}
and one can easily check that the second term involving the
derivative of $F^{*}$ with respect to $\tet$ in expression
\eqref{eq:pb} for $\{F^{*},\LL^{*}\}$ is small compared with the first one
(see Lemma tal in \cite{DelshamsH09}). So, the dominant term in $\eps$ of $\{F^{*},\LL^{*}\}(I, \tet)$ is
\begin{equation}\label{eq:term}
-\left ( I - \frac{l_0}{k_0} \right ) \frac{\partial
\LL^{*}}{\partial \tet} (I, \tet),
\end{equation}
for any $I \in [-1/2,I^{*}_{+}]$.

For Hamiltonian \eqref{eq:example}, by Proposition \ref{prop:prop}, we know
that the reduced Poincar\'e function $\LL^{*}$ in \eqref{eq:Ls} is non-constant in the variable $\tet$ and
indeed, for any $I \in [-1/2,I^{*}_{+}]$, the function
\[ \tet \mapsto \frac{\partial {\LL^{*}}}{{\partial
\tet}} (I,\tet)\] vanishes only for $\tet=0,\pi$, and is negative for
values of $\tet$ between $0$ and $\pi$ and positive between $\pi$
and $2 \pi$. 


Hence, in view of equation \eqref{eq:Escat}, condition \eqref{eq:pois} guarantees
that the scattering map (outer dynamics) moves the invariant tori of the inner
dynamics. It remains to check now that the image under the scattering map of each of these tori intersects transversally another
torus of the foliation. This provides the existence of a transverse heteroclinic connection among these tori. 

Since the different types of tori that appear in our problem
have different quantitative properties and also different expressions, we consider two different cases to check the tranversality condition:
the flat tori region and the big gaps region, introduced in section \ref{sec:inner}.



In the case of flat tori, which by equation \eqref{eq:Fordern} are given approximately by $I=\cte$,
it is clear from expression \eqref{eq:Escat} and \eqref{eq:Lflat}
that the flat invariant tori are mapped under the scattering map to
\[ I - \eps \frac{\partial \LL^{*}}{\partial \tet} (I,\tet) = \cte.\]
Because of the cosinus-like behavior of $\LL^{*}$ described in Proposition \ref{prop:prop} (see also Remark \ref{rem:rpf0}), the image
of this torus intersects transversally other invariant tori with
$I=\cte$, that we know from the previous section that they are at a
distance smaller than $\eps$.


Notice that depending on the point on the torus that we choose to apply the scattering map, we can ``jump'' to different KAM invariant tori with either a higher or a 
lower value of the action $I$. Hence, it is possible to construct several diffusing orbits with increasing or decreasing $I$.


In the big gaps region the computation is more complex because in this case we have primary and secondary KAM tori. Moreover, the invariant KAM tori
are bent near the separatrix. In a connected component of the big gaps regions, around a resonance $I=l_0/k_0$, invariant KAM tori are given approximately by the 
implicit equation
\begin{equation}\label{eq:fe}
F^{*}(I,\tet) = E, 
\end{equation}
for $F^{*}$ as in \eqref{eq:Fsim} and $E$ taking discrete values around $0$. 

When the value of $E$ is positive, equation $F^{*}(I,\tet)=E$ provides two primary invariant tori $\mathcal{T}_{E}^{\pm}$, whereas for $E<0$, it provides a secondary torus 
$\mathcal{T}_E$.  Equation \eqref{eq:fe} defines two smooth surfaces given as graphs of the action $I$ over the angle $\tet$ defined in a certain range:
\begin{equation}\label{eq:graph}
I=f^{*}_{\pm}(\tet,\eps)=\frac{l_0}{k_0} \pm \sqrt{2 (E-\eps a_{k_0,l_0}(\cos(k_0 \tet)-1))}, 
\end{equation}
corresponding to the $+$ and $-$ sign. When $E>0$ these smooth surfaces correspond to the two primary KAM tori $\mathcal{T}_{E}^{\pm}$ and the graph \eqref{eq:graph} is defined 
in the whole domain $[0,2\pi)$. For $E<0$, they correspond to the two components of the secondary KAM torus $\mathcal{T}_{E}$ and the graph \eqref{eq:graph} 
is defined in a domain strictly contained in $(0,2\pi)$.

By expressions \eqref{eq:Escat}, \eqref{eq:term} and \eqref{eq:graph} we have that the
image under the scattering map of a torus satisfying $\{ F(I,\vp,s)=E\}$ is given in first order by the set of points $(I,\tet)$ satisfying
\[F^{*}(I,\tet) \mp \eps \sqrt{2(E - \eps a_{k_0,l_0} (\cos (k_0 \tet) - 1)} \frac{\partial \LL^{*}}{\partial \tet}
\left ( I, \tet \right ) = E.\]


\begin{figure}[ht]
 \begin{center}
 \includegraphics[width=120mm]{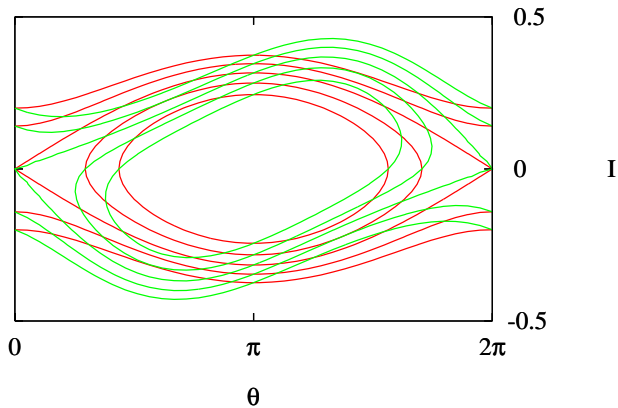}
 \caption{Invariant tori (primary and secondary) in the resonant region around $I=0$ (red
 curves) given implicitly by the level sets of the function
 $F^{*}(I,\tet)$ defined in \eqref{eq:Fsim} with
 $k_0=1$, $l_0=0$ and $a_{1,0}=1/2$. Images of these
 invariant tori (green curves) under the scattering map \eqref{eq:scatt} generated by
 the reduced Poincar\'e function $\LL^{*}(I,\tet)$ given in \eqref{eq:Ls}.}
 \label{fig:biggaps}
 \end{center}
\end{figure}

From this expression, it is clear that if the second term on the left-hand side, namely
\begin{equation}\label{eq:mterm}
\eps\mathcal{M}(\tet,\eps):=\eps \sqrt{2(E - \eps a_{k_0,l_0} (\cos (k_0 \tet) - 1)} \frac{\partial \LL^{*}}{\partial \tet}(I,\tet),
\end{equation}
is non-constant, the image of the invariant torus intersects transversally other invariant 
tori at a distance smaller than $\eps$.

In the case of primary tori (when they are defined as a graph of $I$ over $\tet$ in the whole domain $[0,2\pi)$),  there is an easy way to check this condition. 
Recall first that, by Proposition \ref{prop:prop}, the function $\tet \mapsto \LL^{*}(I,\tet)$ is ``cosinus-like'' and its dominant term possesses two
non-degenerate critical points at $\tet=0,\pi$. Therefore, one can see that the points on the torus
corresponding to $\tet=0, \pi$ remain invariant for the scattering
map, so they are not good for diffusion. However, all the other points
are moved by the scattering map. Indeed, for $\tet \in (0,\pi)$, the
scattering map decreases the value of the action $I$, whereas for
$\tet \in (\pi,2 \pi)$ the scattering map increases it. See Figure \ref{fig:biggaps}. Again, as in the case of the flat tori region, depending on the point on the torus that
the scattering map is applied one can diffuse either increasing or decreasing the value of the action $I$. In this paper we are concerned with diffusion with increasing $I$.

Thus, since
\[ \frac{\partial^2 \LL^{*}}{\partial \tet^2} \neq 0\]
in $\tet=0,\pi$, because these points correspond to non-degenerate extrema of $\tet \mapsto \LL^{*}(I,\tet)$, it is immediate that
\[ \frac{\partial}{\partial \tet} \mathcal{M}(\tet,\eps) \neq 0,\]
in $\tet=0,\pi$. Therefore, there exists an open
neighborhood around these points where the intersection of the invariant tori with its image under the scattering map is transversal.
This is also true for secondary tori in $\tet=\pi$ when $k_0$ is odd. In the other cases and also in general, it is easy to see
in an analogous way that the intersections are transversal
because the ``cosinus-like'' behavior of the function $\tet \mapsto \LL^{*}(I,\tet)$ described in Proposition \ref{prop:prop} guarantees that
\[ \frac{\partial}{\partial \tet} \mathcal{M}(I,\tet) \not \equiv 0.
\]

\begin{figure}
 \begin{center}
 \includegraphics[width=100mm]{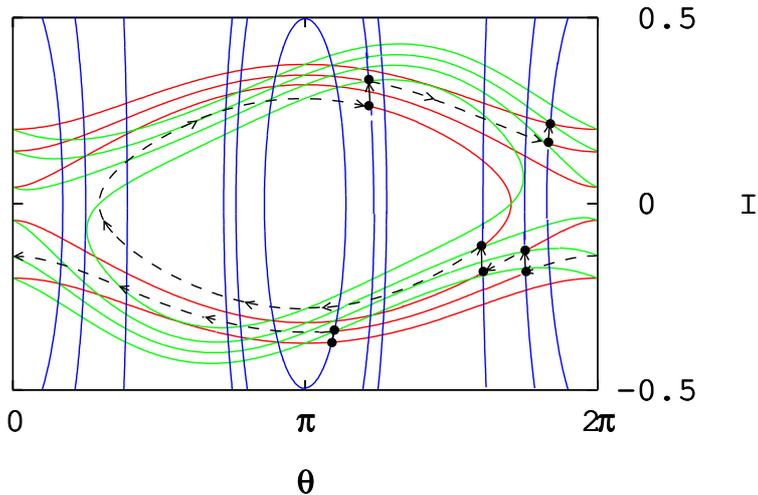}
\caption{Illustration of how to combine the two dynamics to cross the big gaps region. Invariant tori for the inner dynamics (red curves) and invariant sets for the 
outer dynamics (blue curves). Inner dynamics is represented by dashed lines whereas outer dynamics is 
represented by solid lines.}
 \label{fig:twodyn}
 \end{center}
\end{figure}

Up to this point we have shown that it is possible to construct a sequence of invariant tori in the inner dynamics having transverse heteroclinic connections among them, by means 
of the combination of two dynamics. See Figure \ref{fig:twodyn} for an illustration of the combination of the two dynamics. Thus, we have constructed a
transition chain. Finally, using a standard obstruction property (see \cite{DelshamsLS06,FontichM01}) one can show that there exists an orbit that shadows
this transition chain, and Theorem~\ref{teo:main} follows.

\subsubsection*{Acknowledgements}
We are very grateful to R. de la Llave, P. Rold\'an and T. M. Seara
for valuable comments and suggestions.

G. H. has also been supported by the i-Math fellowship ``Contratos Flechados i-Math'', while 
the final version of this manuscript was written.



\end{document}